\documentclass[11pt]{amsart}

\usepackage{amssymb,amscd,fullpage}

\usepackage[all]{xy}
\CompileMatrices

\title{What do DG-categories form?}

\author{Dmitry Tamarkin}

\address{}

\newcommand{\complexes}{{\bf complexes}}

\newcommand{\A}{\mathcal{A}}
\newcommand{\C}{\mathcal{C}}
\newcommand{\F}{\mathcal{F}}

\newcommand{\I}{\mathcal{I}}
\newcommand{\K}{\mathcal{K}}
\renewcommand{\O}{\mathcal{O}}
\renewcommand{\S}{\mathcal{S}}
\newcommand{\T}{\mathcal{T}}
\newcommand{\U}{\mathcal{U}}
\newcommand{\V}{\mathcal{V}}
\newcommand{\W}{\mathcal{W}}
\newcommand{\D}{\mathcal{D}}

\newcommand{\vC}{{\vec{\C}}}
\newcommand{\vc}{{\vec{c}}}
\newcommand{\vF}{{\vec{\F}}}
\newcommand{\vf}{{\vec{f}}}
\newcommand{\vg}{{\vec{g}}}

\newcommand{\vR}{{\vec{R}}}
\newcommand{\vr}{{\vec{r}}}
\newcommand{\vI}{{\vec{I}}}
\newcommand{\vi}{{\vec{\imath}}}
\newcommand{\arI}{\vec{I}}
\newcommand{\vJ}{{\vec{J}}}

\DeclareMathOperator{\shom}{\mathbf{hom}}
\renewcommand{\k}{k}
\newcommand{\vs}{{\sigma}}

\newcommand{\G}{\mathbf{globe}}
\newcommand{\Id}{\mathrm{Id}}
\newcommand{\nat}{\mathbb{N}}

\newcommand{\seq}{{\bf seq}}
\newcommand{\op}{{\text{op}}}
\newcommand{\sets}{{\bf sets}}
\newcommand{\ob}{{\text{\bf Ob}}}
\newcommand{\rhom}{{\text{\bf Rhom}}}
\newcommand{\Rhom}{\rhom}
\newcommand{\triv}{{\text{\bf triv}}}
\newcommand{\Ob}{{\ob}}
\newcommand{\full}{{\text{\bf full}}}
\renewcommand{\bigsqcup}{\coprod}
\newcommand{\chrom}{{\bf CHROM}}
\newcommand{\hoch}{{\bf Hoch}}
\newcommand{\End}{\text{End}}
\newcommand{\E}{\mathcal{E}}

\newtheorem{theorem}{Theorem}[section]

\newtheorem{Lemma}[theorem]{Lemma}
\newtheorem{Corollary}[theorem]{Corollary}
\newtheorem{Proposition}[theorem]{Proposition}
\newtheorem{Theorem}[theorem]{Theorem}

\newtheorem{Definition-Proposition}[theorem]
{Definition-Proposition}

\begin{document}

\begin{abstract}
  We introduce a homotopy 2-category structure on the category of
  2-categories.
\end{abstract}

\maketitle

\section{Introduction}

It is well known that categories form a
2-category: 1-arrows are functors and
2-arrows are their natural transformations.

In a similar way, dg-categories also form a
2-category: 1-arrows $A\to B$ are
dg-functors; given a pair of dg-functors
$F,G:A\to B$ one can define a complex of
their natural transformations $\hom(F,G)$
which naturally generalizes the notion of a
natural transformation  in the usual
setting. Thus, we can use $\hom(F,G)$ as the
space of 2-arrows.

However, this construction has a serious
drawback: the spaces $\hom(F,G)$ are not
homotopically invariant in any way. For
example: let $W:B\to C$ be a  weak
equivalence of dg-categories, we then have a
natural map
$$
\hom(F,G)\to \hom(WF,WG)
$$
which, in general, is not a
quasi-isomorphism of complexes.

Drinfeld \cite{D} proposes another
construction, in which the role of
dg-functors $A\to B$ is played by
$A^\op\times B$-bi-modules. By choosing an
appropriate class of such bi-modules, one
can achieve a good homotopy behavior.
Unfortunately, this class does not contain
identity functors $A\to A$ but only their
resolutions which satisfy the properties of
identity only up to  homotopies.

The goal of this paper is to provide for a
homotopy invariant structure on the category
of dg-category which, on one hand, would be
as close to the 2-category structure as
possible; on the other hand,  this structure
should be free of the above mentioned
drawbacks.

In order to achieve a homotopy invariant
behavior, one has to pass to a derived
version of the notion of a natural
transformation  between two functors. There
is a standard way to do it, in which
$\hom(F,G)$ gets replaced with a certain
co-simplicial complex $ \shom^\bullet(F,G) $
(see Sec. \ref{cosim}).

 As it is
common in such situations, these derived
transformations of functors cannot be
composed as nicely as the usual
transformations of functors do...so, they
don't form a 2-category.

Our result is that, informally speaking,
derived natural transformations form a
certain homotopy version of 2-category. Let
us now sketch the idea of the notion of a
homotopy 2-category, precise definitions
will be given in the main body of the paper.

A good starting point is  to reformulate the
definition of a dg-2-category as follows:

a  small 2-category is

1) a set of objects $C$,  a set of 1-arrows
$\hom(X,Y)$ for every pair $X,Y\in \C$.
These data should form a usual category;

2) a complex of 2-arrows $\hom(F,G)$ for all
one 1-arrows $F,G\in\hom(A,B)$. These data
should have the following structure:

3) given objects $A_0,A_1,\ldots, A_n$ and
1-arrows $F_{ij}:A_i\to A_{i+1}$,
$i=0,1,\ldots,n-1$; $j=0,1,\ldots, m_i$, one
should have a  composition map
\begin{equation}\label{compc}
c:\bigotimes\limits_{ij}\hom(F_{ij};F_{i,j+1})\to
\hom(F,G),
\end{equation}
where $F,G:A_0\to A_n$;
$$F=F_{n-1,0}F_{n-2,0}\cdots F_{10}F_00;$$
$$
G=G_{n-1,m_{n-1}}G_{n-1,m_{n-2}}\cdots
G_{1m_1}G_{0m_0}.
$$

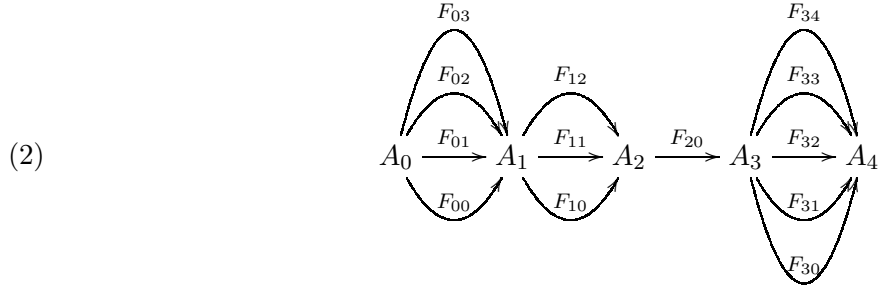
\begin{equation}\label{picc}
\xymatrix{A_0
\ar@/^4pc/[r]^{F_{03}}\ar@/^2pc/[r]^{F_{02}}\ar[r]^{F_{01}}\ar@/_2pc/[r]^{F_{00}}
&
A_1\ar@/^2pc/[r]^{F_{12}}\ar[r]^{F_{11}}\ar@/_2pc/[r]^{F_{10}}&A_2
\ar[r]^{F_{20}}
&A_3\ar@/^4pc/[r]^{F_{34}}\ar@/^2pc/[r]^{F_{33}}
\ar[r]^{F_{32}}\ar@/_2pc/[r]^{F_{31}}\ar@/_4pc/[r]^{F_{30}}&A_4\\}
\end{equation}

There should be a certain coherence axiom
saying that these compositions are closed
under iterations. Instead of giving a
precise formulation, let us consider an
example, as on the picture (\ref{picc}). Let
us split this picture into
 four sub-pictures as follows:

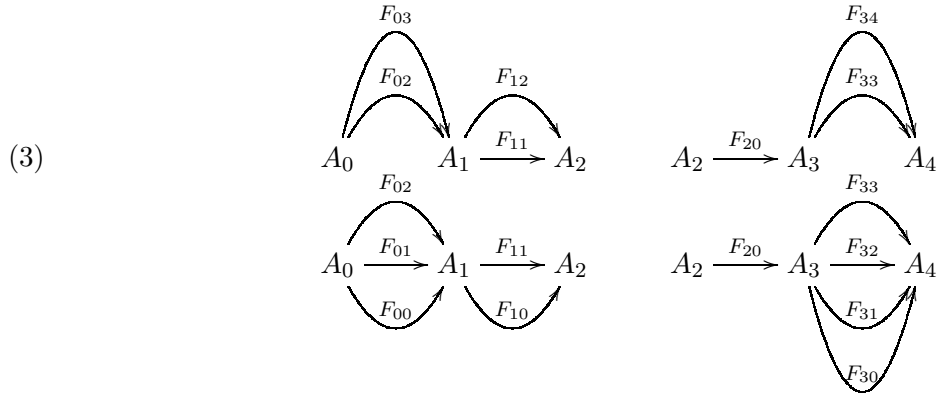
\begin{equation}\label{subpic} \xymatrix{A_0
\ar@/^4pc/[r]^{F_{03}}\ar@/^2pc/[r]^{F_{02}}
&
A_1\ar@/^2pc/[r]^{F_{12}}\ar[r]^{F_{11}}&A_2
&A_2 \ar[r]^{F_{20}}
&A_3\ar@/^4pc/[r]^{F_{34}}\ar@/^2pc/[r]^{F_{33}}
&A_4\\
A_0\ar@/^2pc/[r]^{F_{02}}\ar[r]^{F_{01}}\ar@/_2pc/[r]^{F_{00}}
&
A_1\ar[r]^{F_{11}}\ar@/_2pc/[r]^{F_{10}}&A_2&A_2
\ar[r]^{F_{20}} &A_3\ar@/^2pc/[r]^{F_{33}}
\ar[r]^{F_{32}}\ar@/_2pc/[r]^{F_{31}}\ar@/_4pc/[r]^{F_{30}}&A_4}
\end{equation}

These sub-pictures yield the following
composition maps:

$$
\hom(F_{02},F_{03})\otimes
\hom(F_{11},F_{12})\to
\hom(F_{11}F_{02};F_{12}F_{03});
$$
$$
\hom(F_{33},F_{34})\to
\hom(F_{33}F_{20};F_{34}F_{20});
$$
$$
\hom(F_{00},F_{01})\otimes
\hom(F_{01},F_{02})\otimes\hom(F_{10},F_{11})\to
\hom(F_{10}F_{00};F_{11}F_{02});
$$
$$
\hom(F_{30},F_{31})\otimes
\hom(F_{31},F_{32})\otimes\hom(F_{32},F_{33})\to
\hom(F_{30}F_{20};F_{33}F_{20});
$$

These maps can be composed with the
composition map determined by the following
"quotient-picture" :

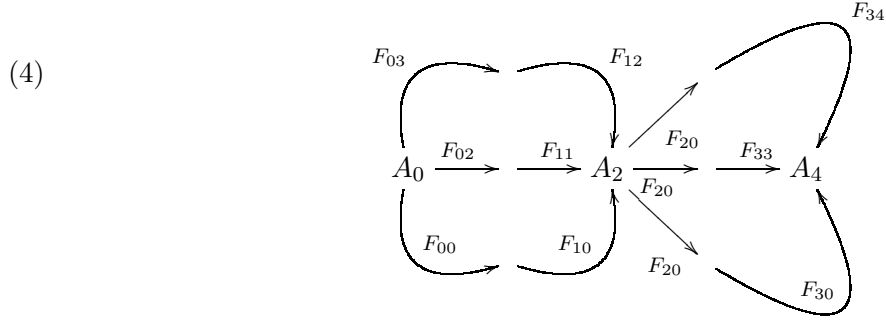
\begin{equation}\label{quotpic} \xymatrix{ &\ar@/^2pc/[rd]^{F_{12}}&&\ar@/^4pc/[rd]^{F_{34}}&\\
A_0
\ar@/^2pc/[ru]^{F_{03}}\ar[r]^{F_{02}}\ar@/_2pc/[rd]^{F_{00}}
& \ar[r]^{F_{11}}&A_2
\ar[ru]_{F_{20}}\ar[r]_{F_{20}}\ar[rd]_>>{F_{20}}
&\ar[r]^{F_{33}}
&A_4\\
&\ar@/_2pc/[ru]^{F_{10}}&&\ar@/_4pc/[ru]^{F_{30}}&}
\end{equation}
\begin{multline*}
\hom(F_{11}F_{02};F_{12}F_{03})\otimes
\hom(F_{10}F_{00};F_{11}F_{02})\otimes\hom(F_{30}F_{20};F_{33}F_{20})
\hom(F_{33}F_{20};F_{34}F_{20})\\
\to
\hom(F_{30}F_{20}F_{10}F_{00};F_{34}F_{20}F_{12}F_{03})
\end{multline*}

so as to get a map
\begin{multline*}\hom(F_{00},F_{01})\otimes\hom(F_{00},F_{01})\otimes\hom(F_{01},F_{02})\otimes\hom(F_{02},F_{03})\\
\otimes\hom(F_{10},F_{11})\otimes\hom(F_{11},F_{12})\\
\otimes\hom(F_{30},F_{31})\otimes\hom(F_{31},F_{32})
\otimes\hom(F_{32},F_{33})\otimes\hom(F_{33},F_{34})\\
\to
\hom(F_{30}F_{20}F_{10}F_{00};F_{34}F_{20}F_{12}F_{03})
\end{multline*}

The coherence axiom then requires that {\em
the map that we have just constructed should
coincide with the map (\ref{compc})
determined by the picture (\ref{picc})}.

 This definition can be homotopized in the
standard way: we define a homotopy
2-category as:

a collection of data 1),2) (same as in the
above definition) with 3) modified as
follows:

h3) for each collection of non-negative
integers $m_0,m_1,m_2,\ldots,m_n$ there
should be given a contractible complex
$\O(m_0,m_1,\ldots,m_n)$ and a map

$$
c:\O(m_0,m_1,\ldots,m_n)\otimes\bigotimes\limits_{ij}\hom(F_{ij};F_{i,j+1})\to
\hom(F,G),
$$
where the notations are the same as in 3).

In order to formulate the coherence axiom we
need some operad-like structure on the
collection of complexes
$\O(m_0,m_1,\ldots,m_n)$. Let us briefly
discuss this structure.

First, with every picture $P$ as in
(\ref{picc}), one naturally associates a
complex $\O(P)$ (example: for the picture
from (\ref{picc}), $\O(P):=\O(3,2,0,4)$,
where the numbers represent the numbers of
arrows in each column of the picture); next,
suppose we have a subdivision of a picture
$P$ into a number of sub-pictures
$P_1,P_2,\ldots,P_k$ with the corresponding
quotient-picture $Q$ (we do not define the
precise meaning of these words hoping that
the spirit can be felt from the above
example of a subdivision of the picture
(\ref{picc}) into four sub-pictures
(\ref{subpic}) with the corresponding
quotient-picture (\ref{quotpic})).

We then should have a composition map
\begin{equation}\label{twoop}
\O(P)\otimes \O(P_1)\otimes \O(P_2)\cdots
\O(P_k)\to \O(Q)
\end{equation}

Having these maps, one can formulate the
coherence axiom in this new setting in a
natural way.

The maps (\ref{twoop})  along with certain
natural associativity properties constitute
a so-called structure of 2-operad
\cite{Bat}; we  will reproduce a precise
definition.

Our main result is that {\em dg-categories
form a homotopy 2-category} in which objects
are dg-categories, 1-arrows are functors and
the complexes of 2-arrows are defined using
the derived version of the complex of
natural transformations.

We conclude the paper with an observation
that this result immediately implies that
for every category $\C$, the complex
$\rhom(\Id_C,\Id_C)$ (the homotopy center of
$\C$) is an algebra over the above mentioned
2-operad $\O$. A result from \cite{Bat}
implies that {\em an algebra structure over
any contractible 2-operad ($\O$ is such)
implies a structure of an algebra over some
resolution of the chain operad of little
disks}, thus yielding another proof of
Deligne's conjecture on Hochschild cochains
\cite{MS},\cite{KS},\cite{V}.

The plan of the paper is as follows. We
begin with defining a co-simplicial complex
of natural transformations
$\shom^\bullet(F,G)$ for every pair of
dg-functors $F,G:A\to B$. By taking the
realization, one gets a complex
$\rhom(F,G):=|\shom^\bullet(F,G)|$ which we
use as a replacement for the naive complex
$\hom(F,G)$.

 Next, we introduce
some combinatorics in order to describe
pictures like the one in (\ref{picc}).

Next, we make definitions of a 2-operad
(which is equivalent to that in \cite{Bat})
and a homotopy 2-category.

After that, we proceed to constructing a
homotopy 2-category of dg-categories.  It
turns out to be more convenient to start
with constructing a certain structure on the
co-simplicial complexes
$\shom^\bullet(F,G)$, without passing to the
realization. This structure will be given in
terms of a collection of certain
poly-simplicial sets so that one  can study
them using some topology.

Finally, using the realization functors,
this structure will be converted to a
homotopy 2-category structure in  which the
complexes of 2-arrows are $\rhom(F,G)$.

We conclude by showing that this result
coupled with Batanin's theorem on
contractible 2 operads readily implies
Deligne's conjecture.

I would like to thank  A. Beilinson for
statement of the problem and  M. Batanin for
explaining me his results.

\section{Conventions, notations}
\subsubsection{Ordered sets} Any finite
non-empty totally ordered set will be called
{\em an ordinal}.

Given a non-negative integer $n$, we denote
by $[n]$ the ordinal $\{0<1<\cdots<n\}$.

Given an ordinal $I$ we denote by $m_I$ its
minimum and by $M_I$ its maximum.

Denote by $\vI$ the set of all pairs
$\vec{\imath\jmath}$, where $i,j\in I$ and $j$ is the
immediate successor of $i$. We have an
induced total order on $\vI$, but $\vI$ may
be empty.  We have natural projections
$s,t:\vI\to I$; $s(\vec{\imath\jmath})=i$;
$t(\vec{\imath\jmath})=j$.

Given  $a,b\in I$, $a\leq b$, we define the
interval $[a,b]\subset I$ in the usual way.

A monotonous(:=non-decreasing) map of
ordinals  $f:I\to J$ will be called {\em
dominant} if  $f(m_I)=m_J$; $f_(M_I)=M_J$.

\section{Functors between dg categories,
and their natural transformations}
\subsubsection{} Let $A$ be a dg-category. Let
$X:I\to A$ be a family of objects in $A$
indexed by an ordinal $I$.

Set $A(X):=\bigotimes\limits_{\vi\in \arI}
A(s(\vi);t(\vi))$. If $I$ is a one-element
ordinal, we set $A(X)=\k$.

Given a dominant monotonous map $k:J\to I$,
we have a natural map
$$
k^*:A(X)\to A(X\circ k)
$$
\subsubsection{Definition of the co-simplicial complex of natural transformations}\label{cosim} Let $A$, $B$
be small dg categories,
$$
F,G:A\to B
$$
functors.

Let $I$ be a finite non-empty totally
ordered set. Set
$$
\shom^I(A,B):=\prod\limits_{X:I\to A}
\hom_\k\Big(A(X);\hom_B\big(F(X(m_I));G(X(M_I))\big)\Big).
$$

Let $I'$ be obtained from $I$ by adding two
elements $m',M'$ such that $m'<I<M'$.  Let
$X':I'\to A$, let $X$ be the restriction of
$X'$ onto $I$ We then have
\begin{equation}\label{id0}
A(X')\cong A(X(M_I);X'(M'))\otimes
A(X)\otimes A(X'(m');X(m_I)) \end{equation}

We  have a natural map
$$
\shom^I(F,G)\to \shom^{I'}(F,G)
$$

such that the chain $\Phi\in \shom^I(F,G)$
is mapped into a chain $\Phi'$ according to
the rule
$$
\Phi'(\omega\otimes U\otimes
\alpha)=G(\omega)\circ \Phi(U)\circ
F(\alpha),
$$
where $\omega\in A(X(M),X'(M'))$; $U\in
A(X)$; $\alpha \in A(X'(m'),X(m))$, and we
use the identification (\ref{id0}).
\subsubsection{cosimplicial structure} Let
$\Delta$ be the category of ordinals and
their monotonous (=non-decreasing) maps. We
are going to endow the collection of spaces
$\shom^I(F,G)$ with a structure of  functor
$$
\Delta\to \complexes.
$$

 Let
$\vs:I\to J$ be a monotonous map. Define a
map
$$
\vs_*:\shom^I(F,G)\to \shom^J(F,G)
$$
as follows.

Let $\vs':I'\to J$ be the extension of $\vs$
which sends $m'$ and $M'$ to the minimum and
the maximum of $J$ respectively.

Set $\vs_*\Phi(U):=\Phi'((\vs')^*U)$, where
$U\in A(X)$ for some $X:J\to A$.

This way we get the desired co-simplicial
structure.

\subsection{ Total complex of a cosimplicial comlplex}
\label{simplchain} It is well known that
given a co-simplicial complex one can
produce its total complex by applying
alternated sums of co-face maps.

We will use  a slightly different definition
of this total complex.

Let $I$ be an ordinal, let $\Delta^I$ be the
simplex whose vertices are labelled by $I$.

Let $C_*(\Delta^I)$ be its reduced chain complex. Let
$\S^{*}(I):=C_{-*}(\Delta^I)$. It is clear that $ \S^*(\bullet) $ is a
co-simplicial complex (here $\bullet$ stands for the "co-simplicial"
argument).

We will denote this co-simplicial complex simply by $\S$.

Given a co-simplicial complex $K$ we can form a complex
$\hom_{\Delta}(\S,K)$ which will be also denoted by $|K|$.

Thus, we can construct a complex
$|\shom^\bullet(F,G)|$ which will  be
denoted by $\rhom(F,G)$.

\subsubsection{}
\label{consthom} We also have a "naive"
notion of the complex of natural
transformations of two functors.

Indeed, given a pair of functors $F,G:A\to
B$, define the complex
$$
\hom(F,G)
$$

as the equalizer of the diagram
$$\xymatrix{
\hom(F,G)\ar[r] &
\shom^{[0]}(F,G)\ar@<1ex>[r]^{d_1}\ar@<-1ex>[r]_{d_0}
& \shom^{[1]}(F,G)}
$$
where $d_0,d_1$ are the co-face maps.

We can define a constant co-simplicial
complex
$$
\hom^\bullet(F,G),
$$
where
$$
\hom^I(F,G):=\hom(F,G).
$$

We then have a  co-simplicial map
$$
\hom^\bullet(F,G)\to \shom^\bullet(F,G).
$$
\section{Some combinatorics} We want
to find an algebraic structure naturally
possessed by complexes $\rhom(F,G)$.

This structure will be given in terms of a
family of poly-linear maps between these
complexes and some relations between them.

In order to formulate this structure we need
some combinatorics.

\subsection{Combinatorial data}
\subsubsection{2-ordinals, 2-trees}
By definition, a 2-ordinal $\U$ is a
collection of the following data:

--- an ordinal $\C_\U$;

--- for each $\vc\in \vC_\U$, an ordinal
$\F_{\U,\vc}$.

2-ordinals are meant to represent pictures
of the type shown below:\bigskip

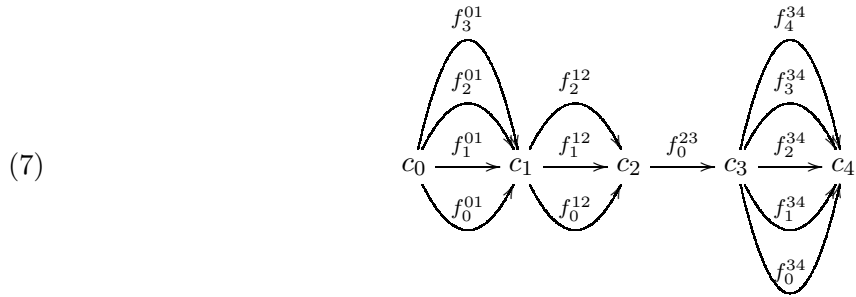
\begin{equation}\label{pic} \xymatrix{c_0
\ar@/^4pc/[r]^{f^{01}_3}\ar@/^2pc/[r]^{f^{01}_2}\ar[r]^{f^{01}_1}\ar@/_2pc/[r]^{f^{01}_0}
&
c_1\ar@/^2pc/[r]^{f^{12}_2}\ar[r]^{f^{12}_1}\ar@/_2pc/[r]^{f^{12}_0}&c_2
\ar[r]^{f^{23}_0}
&c_3\ar@/^4pc/[r]^{f^{34}_4}\ar@/^2pc/[r]^{f^{34}_3}
\ar[r]^{f^{34}_2}\ar@/_2pc/[r]^{f^{34}_1}\ar@/_4pc/[r]^{f^{34}_0}&c_4\\}
\end{equation}

\bigskip

This picture corresponds to the following
2-ordinal:

--- $\C=\{c_0<c_1<c_2<c_3<c_4\}$;

---
$\F_{\vec{c_0c_1}}=\{f^{01}_0<f^{01}_1<f^{01}_2<f^{01}_3\}$;

---
$\F_{\vec{c_1c_2}}=\{f^{12}_0<f^{12}_1<f^{12}_2\}$;

---
$\F_{\vec{c_2c_3}}=\{f^{23}_0\}$;

--
$\F_{\vec{c_3c_4}}=\{f^{34}_0<f^{34}_1<\cdots<f^{34}_4\}
$

A picture of such a type can be drawn for
any 2-ordinal in the obvious way.

 We denote
$$\vF_\U:=\sqcup_{\vc\in \vC}
\vF_{\U,\vc}.$$

This set is in 1-to-1 correspondence with
the set of 2-cells of the picture
corresponding to $\U$.

We then have an obvious monotonous  map of
finite totally ordered sets.
$$
\pi:\vF_{\U}\to \vC_\U.
$$

According to \cite{Bat}, let us call any map
of finite totally ordered sets {\em a 2
stage tree} or, shortly, {\em a 2-tree}.

We will denote  2-trees by one letter, say
$t$. We will refer to the elements of $t$
as:
$$
\pi_t:\vF_t\to \vC_t
$$

We have shown how, given a 2-ordinal $\U$,
one constructs a 2-tree $\pi_U:\vF_\U\to
\vC_\U$. Denote this 2-tree by $t_\U$, so
that:
$$
\pi_{t_\U}:=\pi_\U;\ \vF_{t_\U}:=\vF_\U;
$$
$$
\vC_{t_\U}:=\vC_\U.
$$

It is clear that a 2-ordinal is defined
up-to a canonical isomorphism by its 2-tree.

\subsubsection{} Given a 2-ordinal $\U$ we
can construct a strict 2-category $[\U]$,
the universal one among the 2-categories $V$
possessing the following properties:

--- $\ob V= \C$;

--- there are fixed  maps
$\F_{\vec{c_1c_2}}\to \Ob\hom_V(c_1,c_2)$
for all $\vec{c_1c_2}\in \vC$. For an $f\in
\F_{\vec{c_1c_2}}$ we denote by the same
symbol the corresponding object in
$\hom_V(c_1,c_2).$

--- For each $\vec{f_1f_2}\in
\vF_{\vec{c_1c_2}}$ we have a fixed element
in $\hom_{\hom_V(c_1,c_2)}(f_1,f_2)$.

\bigskip

This 2-category has a clear meaning in terms
of the picture (\ref{pic})

Objects are $c_0,c_1,\cdots$; the space of
maps $c_i\to c_j$  is non-empty iff $c_i\leq
c_j$, in which case an arrow $c_i\to c_j$ is
just a directed path from $c_i$ to $c_j$;
let us define a partial order on  the space
of such paths by declaring that one  path is
less or equal to another iff the former lyes
below the latter. We then have a category
structure on $\hom(c_i,c_j)$ produced by the
just defined poset of paths (each arrow goes
from a smaller object to a greater one).

Here is a more formal
 description.
Given $c,c'\in \C$, we have

 1) $\hom(c,c')=\emptyset$ if $c>c'$;

 2) $\hom(c,c)=\{\Id_c\}$;

 3) if $c<c'$, then
 $$\Ob\hom(c,c'):=\prod\limits_{\vc\in
 \vec{[c,c']}} \F_{\vc}.
$$

Given $f^1,f^2\in \hom(c,c')$;

$$
f^k=\{f^k_{\vc}\}_{\vc\in \vec{[c_1c_2]}},\
k=1,2;
$$

we have a unique arrow $f_1\to f_2$ iff
$$
f^1_\vc\leq f^2_\vc
$$
for all $\vc\in \vec{[c_1c_2]}$. Thus, the
set $\hom(c_1,c_2)$ is partially ordered and
has the least and the greatest element.

\subsubsection{} We will often need a
special 2-ordinal called {\em globe} and
denoted by $\G$. We define $\G$ by
$$
\C_\G=\{c_0<c_1\};
$$
$$
\F_{\vec{c_0c_1}}=\{f_0<f_1\}
$$
\bigskip
$$
\xymatrix{\bullet\ar@/^1pc/[r]\ar@/_1pc/[r]&\bullet\\}
$$
\bigskip

Any 2-ordinal isomorphic to $\G$ will also
be called a globe.

\subsubsection{Balls in a 2-ordinal}
 Given a 2-ordinal $\U$ define {\em a
ball} in $\U$ as any 2-ordinal $\U'$
 the form:

 --- $\C_{\U'}$ is an interval in
 $\C_\U$;

 --- for each $\vc\in \vC_{\U'}$,
 $\F_{\U',\vc}$ is an interval in
 $\F_{\U,\vc}$.

 We then see that
$$
[\U']\subset [\U]
$$
is a full subcategory.

The set of all balls in $\U$ is partially
ordered; each minimal ball is a globe; the
set of all these minimal balls is naturally
identified with $\vF_\U$.

\subsubsection{} We write $[\U]_1$ for the
underlying usual category of $\U$.

\subsubsection{Maps of 2-ordinals} Let
$\U,\V$ be 2-ordinals. {\em A map} $P:\U\to
\V$ is a 2-functor $[P]:[\V]\to [\U]$
satisfying:

1) the induced map $[P]:\C_\V\to \C_\U$ is
dominant (= monotonous and preserves the
minima and the maxima);

2) for all $c_1<c_2$, $c_1,c_2\in \C_V$ the
induced map
$$
\hom_{[\V]}(c_1,c_2)\to
\hom_{[\U]}(P(c_1),P(c_2))
$$

preserves the least and the greatest
elements.

With this definition of a map, 2-ordinals
form a category.

Any globe is a terminal object in this
category.

\subsubsection{Inverse images of balls}
Let $P:\U\to \V$ be a map of 2-ordinals and
$\V'\subset \V$ a ball in $\V$. Define
$P^{-1}\V'=:\U'$  as a unique ball in $\U$
satisfying:

 there exists a map of
2-ordinals $P':\U'\to \V'$ fitting into a
commutative diagram:
$$
\xymatrix{ [\U]&[\V]\ar[l]^{[P]}  \\
            [\U']\ar[u]&[\V']\ar[l]^{[P']}\ar[u]
            }
$$
with the vertical arrows being the natural
inclusions.

 Consider a pictorial example:

Let $\U$ be the same ordinal as in
(\ref{pic}):

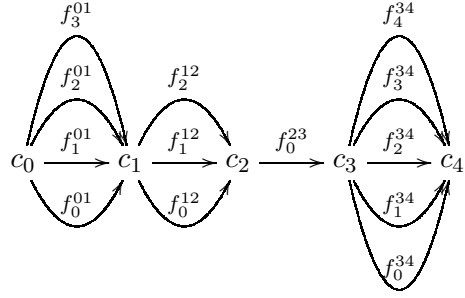
\bigskip

$$ \xymatrix{c_0
\ar@/^4pc/[r]^{f^{01}_3}\ar@/^2pc/[r]^{f^{01}_2}\ar[r]^{f^{01}_1}\ar@/_2pc/[r]^{f^{01}_0}
&
c_1\ar@/^2pc/[r]^{f^{12}_2}\ar[r]^{f^{12}_1}\ar@/_2pc/[r]^{f^{12}_0}&c_2
\ar[r]^{f^{23}_0}
&c_3\ar@/^4pc/[r]^{f^{34}_4}\ar@/^2pc/[r]^{f^{34}_3}
\ar[r]^{f^{34}_2}\ar@/_2pc/[r]^{f^{34}_1}\ar@/_4pc/[r]^{f^{34}_0}&c_4\\}
$$
\bigskip

and let $\V$ be defined by the following
picture
\bigskip
$$
 
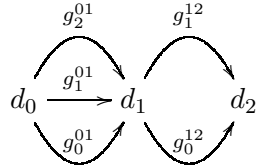
\xymatrix{d_0
\ar@/^2pc/[r]^{g^{01}_2}\ar[r]^{g^{01}_1}\ar@/_2pc/[r]^{g^{01}_0}
&
d_1\ar@/^2pc/[r]^{g^{12}_1}\ar@/_2pc/[r]^{g^{12}_0}&d_2}
$$

Consider a map $P:\U\to \V$ such that the
corresponding map $[P]:[\V]\to [\U]$ is as
follows:

$$
[P]d_0=c_0;\ [P]d_1=c_2;\ [P]d_2=c_4;
$$
$$
[P]g_0^{01}=f^{12}_0f^{01}_0;\
[P]g_1^{01}=f^{12}_1f^{01}_0\
[P]g_2^{01}=f^{12}_2f^{01}_3;
$$
$$
[P]g_0^{12}=f^{34}_0f^{23}_0;
$$
$$
[P]g_1^{12}=f^{34}_4 f^{23}_0.
$$

Let us label globes in $\V$ by $I,II,III$ as
shown on the picture:

$$
 
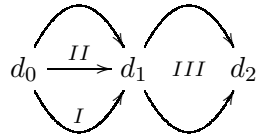
\xymatrix{d_0
\ar@/^2pc/[r]^{}\ar[r]^{II}\ar@/_2pc/[r]^{I}
&
d_1\ar@/^2pc/[r]^{}\ar@/_2pc/[r]\ar@{}[r]|{III}&d_2}
$$

\bigskip

The preimages of these globes are then as
follows

$$ 
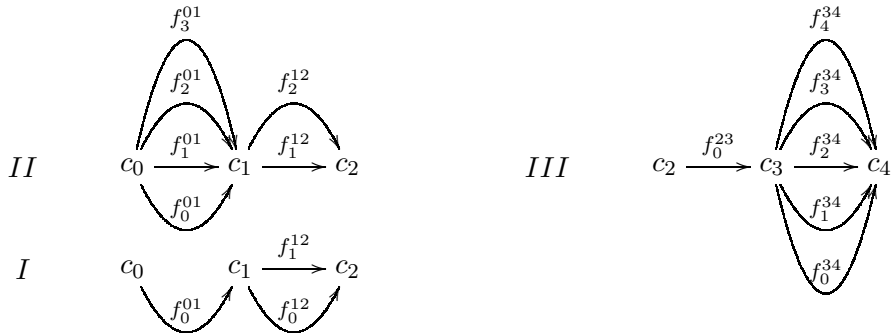
\xymatrix{II&c_0
\ar@/^4pc/[r]^{f^{01}_3}\ar@/^2pc/[r]^{f^{01}_2}\ar[r]^{f^{01}_1}\ar@/_2pc/[r]^{f^{01}_0}
&
c_1\ar@/^2pc/[r]^{f^{12}_2}\ar[r]^{f^{12}_1}&c_2&&III&c_2
\ar[r]^{f^{23}_0} &c_3
\ar@/^4pc/[r]^{f^{34}_4}\ar@/^2pc/[r]^{f^{34}_3}
\ar[r]^{f^{34}_2}\ar@/_2pc/[r]^{f^{34}_1}\ar@/_4pc/[r]^{f^{34}_0}&c_4\\
I&c_0 \ar@/_2pc/[r]^{f^{01}_0} &
c_1\ar[r]^{f^{12}_1}\ar@/_2pc/[r]^{f^{12}_0}&c_2
&&&&&}
$$

\subsubsection{Maps of the corresponding
2-trees} According to Batanin, we define a
map of 2-trees $P:t_1\to t_2$ as a
commutative diagram
$$\xymatrix{ \vF_{t_1}\ar[d]^{\pi_{t_1}}\ar[r]^{P_\vF} & \vF_{t_2}\ar[d]^{\pi_{t_2}}\\
             \vC_{t_1}\ar[r]^{P_\vC} &
             \vC_{t_2}}$$
with all its arrows being monotonous.

\subsubsection{}
A map of ordinals $\U\to \V$ naturally
induces a map of the sets of two-cells: $$
\F_\U\to \F_V;
$$
it is not hard to see that this map lifts to
a map of the corresponding 2-trees. We are
going to give a formal definition of this
map.

 Given
a map of 2-ordinals $P:\U\to \V$, we define
an induced   map of the associated 2-trees:
$$
P^t:t_\U\to t_V
$$
as follows:

1) define the map

$$
P^t_{\vC}:\vC_\U\to \vC_V.
$$

For $\vec{c_1c_2}\in \vC_U$ we set
$P_{\vC}(\vec{c_1c_2})=\vec{d_1d_2}$ iff
$$
[c_1c_2]\subset [[P](d_1),[P](d_2)].
$$

2)  Given a globe $\vf\in \vF_\U$ define its
image $P^t_\vF(\vf)=:\vg$ in $\vF_\V$ as a
unique globe such that  the ball $P^{-1}\vg$
contains $\vf$.
\subsubsection{} We see that this way we
have constructed a category of 2-ordinals
and a category of 2-trees and a functor
between them; it is easy to see that this
functor is an equivalence.

\subsubsection{Diagrams } Given a 2-ordinal
$\U$ and a category $\C$, {\em a
$\U$-diagram} in $\C$ is a functor
$$
\D:[U]_1\to \C
$$

 Given a map of 2-ordinals
$P:\U\to \V$  and a $\U$-diagram $\D$, it
naturally restricts to produce $P^{-1}\vf$-
diagrams $\D|_{P^{-1}\vf}$, where $\vf\in
\vF_\V$), and a $\V$-diagram $P_*\D$.  These
induced diagrams come from the  obvious
functors

$$
[P^{-1}\vf]\hookrightarrow [\U];
$$
$$
[P]:[\V]\to [U].
$$
\section{2-operads and their algebras}
 We are going to adjust notions of an operad
 and an algebra over an operad so that they
 work in our setting. In the usual setting,
 given an operad, we define its action on a
 complex; in our situation, instead of one complex,
 we have a family of complexes: a complex
 $\rhom(F,G)$ for each globe formed by the
 pair of dg-categories $A,B$ and a pair of
 dg-functors $F,G:A\to B$. We abstract this
 situation by introducing a notion of
 a $\C$-complex.  Next,  following \cite{Bat},
 we define the notion of a 2-operad, and, lastly,
 the notion of a structure of an $\O$- algebra
on a $\C$-complex, where $\O$ is a 2-operad.

Using these notions, we will be able to
formulate the definition of a homotopy
2-category as an algebra structure over a
contractible 2-operad.

 \subsection{$\C$-complexes and their tensor products}

We fix a small category in sets $\C$. Let
$\C_0$ be the set of objects in $\C$ and
$\C_1$ be the set of its arrows. Let
$s,t:\C_1\to \C_0$ be the source and target
maps

 We define
{\em a globe in $\C$} as any $\G$-diagram in
$\C$ (= a pair of objects in $\C$ and a pair
of arrows between these objects).

Let $\C_2$ be the set of all globes in $\C$.
We have obvious maps
$$
s',t':\C_2\to \C_1:
$$

given a globe $g$
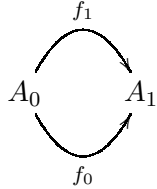
$$\xymatrix{
A_0\ar@/^2pc/[r]^{f_1}\ar@/_2pc/[r]_{f_0} &
A_1}
$$
where $A_0,A_1\in \C$; $f_0,f_1:A_0\to A_1$,
we set $s'(g)=f_0$; $t'(g)=f_1$. It is clear
that $ss'=ts'$; $st'=tt'$, and $\C_2$ is the
universal set with these properties.

\subsubsection{}

We define {\em a $\C$-complex} as a family
of complexes parameterized by $\C_2$.

\subsubsection{Tensor product of
$\C$-complexes} Given a 2-ordinal $\U$ and a
$\vF_\U$-family  of $\C$- complexes
$\K=\{K_\vf\}_{\vf\in \vF_\C}$,   we define
a new $\C$-complex
$$
\bigotimes^\U \K:=\bigotimes_{\vf \in
\vF_\U}^\U K_\vf
$$
as follows.

Pick a globe $g\in C_2$ and define
$$
(\bigotimes^\U \K)_g:=\bigoplus\limits_{\D|
p_*\D=g} \bigotimes_{\vf\in \vF_\U}
K_{(\D|_{\vf})},
$$
where  $p:\U\to \G$ is the terminal map,
$\D$ is any $\U$-diagram in $\C$ with
$p_*\D=g$, and $\vF_\U$ is identified with
the set of globes in $\U$ so that $\D|_\vf$
is a globe in $\C$.
\subsubsection{} Given a map of 2-ordinals
$\U\to \V$, we have a canonical isomorphism
$$
\bigotimes^\V_{\vf\in \vF_\V} \bigotimes
^{P^{-1}\vf}_{\vg\in \vF_{P^{-1}\vf}} K_\vg
\to \bigotimes^\U_{\vf\in \vF_\U} K_\vf;
$$

call this isomorphism {\em a constraint.}

Given  a chain of maps of 2-ordinals $\U\to
\V\to \W$, we get the associativity property
of this constraint.

This can be reformulated so that the
category of $\C$-complexes becomes a
category with two monoidal structure such
that one of them distributes over the other
but we won't need it in this paper.

\subsubsection{} Fix  a set $\chrom$ to be called
the set of colors. Let us also fix
 family  of
$\C$-complexes $\K:=\{K_\chi\}_{\chi\in
\chrom}$.

Let $\U$ be a 2-ordinal. Define {\em an
$\chrom$- coloring $\chi$} of $\U$ as a
prescription of a color $c^\chi_\vf\in
\chrom$ for each $\vf\in \vF_\U$ and an
additional color $c^\chi\in \chrom$.

A 2-ordinal endowed with a coloring will be
called a colored 2-ordinal. Given a colored
2-ordinal $\U':=(\U,\chi)$, we have a
complex
$$
\full_\K(\U'):=\hom([\bigotimes^\U_{\vf\in
\vF_\U} K_{c^\chi_\vf}];K_{c^\chi})
$$

These complexes form an algebraic structure
 called a colored 2-operad. We are going to
 define this notion.

\subsection{Colored 2-operads} We need a
notion of a map of colored 2-ordinals:

Let $\U'=(\U,\chi_\U)$; $\V'=(\V,\chi_\V)$
be colored ordinals. We say that we have a
map $P':\U'\to \V'$ if:

--- we are given a map $P: \U\to \V$  of 2-ordinals;

--- $c^{\chi_\U}=c^{\chi_\V}$.

Given a globe $\vf\in \vF_\V$, we then have
a natural coloring $\chi^\vf$ on
$P^{-1}\vf$: we set
$c^{\chi^\vf}=c^{\chi_{\V}}_\vf$;
$(c^{\chi_\vf})_\vg:=(c^{\chi_\U})_\vg$.

\subsubsection{Definition of a colored
2-operad} A $\chrom$-colored operad $\O$ is
a collection of the following data:

--- a functor $\O$ from the  isomorphism groupoid
of the category of colored 2-ordinals to the
category of complexes;

--- given a map $P:\U'\to \V'$ of colored ordinals,
set
$$
\O(P):=\bigotimes_{\vf\in \vF_\V}
\O(P^{-1}\vf)
$$
 we then should have a map
$$
\circ_P:\O(P)\otimes \O(\V')\to \O(\U').
$$

Axioms:

the first axiom asks for covariance of the
map $\circ_P$ under isomorphisms of
2-ordinals.

 In order to formulate the next axiom, note
that given a chain of maps of 2-ordinals
$$
\U\stackrel P\to \V\stackrel Q\to \W,
$$
we have a natural map
$$\circ(Q,P):\O(P)\otimes \O(Q)\to \O(QP);
$$

indeed, for each $\vf\in \vF_\W$ we have
induced maps
$$
P_\vf:(QP)^{-1}\vf\to Q^{-1}\vf,
$$
and we can define our map as follows
$$
\O(P)\otimes \O(Q)=\bigotimes_{\vf\in
\vF_\W} [\O(P_\vf)\otimes \O(Q^{-1}\vf)]\to
\bigotimes_{\vf\in \vF_\W}
\O((QP)^{-1}\vf)=\O(QP).
$$

The property then reads that the maps
$\circ(Q,P)$ should be associative in the
obvious way.
\subsubsection{} It is immediate that the
complexes $\full_\K(U')$ form a colored
2-operad $\full_\K$

\subsubsection{} Given an $\chrom$-colored 2-operad
$\O$, we define an $\O$-algebra as an
$\chrom$-family $\K$ of $\C$-complexes along
with a map of colored 2-operads
$$
\O\to \full_\K
$$
\subsubsection{} Define a (non-colored)
2-operad as a colored operad with the set of
colors to have only one element.
\subsection{Main Theorem}
 DG-categories and their functors form a
 category. Fix a small sub-category $\C$ in
 this category.

Given a globe $g$ in $\C$, we have defined a
complex $\Rhom(g)$.  These complexes form a
$\C$-complex, to be denoted $\Rhom$.
Likewise we have the functor of usual
homomorphisms $\hom(g)$, these also form a
$\C$-complex $\hom$. We have a natural map
of $\C$-complexes
$$
\hom\to \Rhom.
$$

We know that  the pair $(C,\hom)$ is
naturally a 2-category. This can be
formulated in our language as follows.
Define a trivial 2-operad $\triv$ as
follows:

set $\triv(\U)=\k$ for each 2-ordinal $\U$
and demand that all structure maps preserve
$1\in \k$.

Then the 2-categorical structure on
$\C,\hom$ amounts to the fact that we have a
$\triv$-algebra structure on $\hom$.

\subsubsection{Formulation of a theorem}
Define a notion of {\em a contractible
2-operad} as a collection of the following
data:

--- a 2-operad $\O$;

--- a quasi-isomorphism of 2-operads $\O\to
\triv$.

\begin{Theorem}\label{maintheor} There exists a
contractible 2-operad $\O$ and $\O$-algebra
structures on $\hom$ and $\Rhom$ such that:

1) the map $\hom\to \Rhom$ is a map of
$\O$-algebras;

2) the $\O$-algebra structure on $\hom$ is
the pull-back of the $\triv$-structure via
the structure map $\O\to \triv$.

\end{Theorem}

The rest of the paper is devoted to proving
this theorem.
\section{Constructing a colored 2-operad which
 acts on $\shom^\bullet$}

Let $\nat$ be the set of isomorphism classes
of ordinals; $\nat=\{[0],[1],[2],\ldots\}$.

We then have a $\nat$-family of
$\C$-complexes $I\mapsto \shom^I $, $I\in
\nat$.

We will start with a construction of a
colored 2-operad in the category of sets
naturally acting on $\shom^\bullet$. This
operad will be denoted by $\seq$.

{\em By default, all colorings are
$\nat$-colorings.}

\subsection{Construction of $\seq$}

Let $(\U,\chi)$ be a colored 2-ordinal. Let
$$
\pi:\vF_\U\to \vC_\U$$
be the induced
2-tree.

Let us  use the following notation for the
ordinals which determine the coloring:
$$
I_\vf:=c^\chi_\vf;
$$
$$
J:=c^\chi.
$$

 Define a set $\seq(\U)$ whose each
element is a collection of the following
data:

A) a total order on
$\I:=\I_\U:=\bigsqcup_{\vf\in\vF} I_\vf$;

B) a monotonous  map $W:\I\to J$.

The conditions are:

1) the total order  on $\I$ agrees with the
orders on  each $I_\vf$;

2) if $i_1,i_2\in I_\vf$ and $i_1<i<i_2$,
$i\in I_{\vf_1}$ then $\pi(\vf_1)<\pi(\vf)$;

3) if $\pi(\vf_1)=\pi(\vf_2)=\vc$ and
$\vf_1<\vf_2$ in the sense of the order on
$\vF_\vc$, then $I_{\vf_1}<I_{\vf_2}$ with
respect to the order on $\I$.

\subsubsection{Compositions}  Let $P:\U'\to
\V'$ be a map of $\nat$-colored ordinals.
Define the structure map
$$
\circ_p: \prod_{\vf\in \vF_\V}
\seq(P^{-1}\vf)\times \seq(\V)\to \seq(\U).
$$

Let us pick elements $\lambda_\vf\in
\seq(P^{-1}\vf)$; $\lambda\in \seq(V)$ and
define their composition.

We have
$$
\I_\U=\bigsqcup\limits_{\vf\in \vF_\V}
\I_{P^{-1}\vf}
$$

The elements $\lambda_\phi$ define  total
order on $\I_{P^{-1}\vf}$ and monotonous
maps
$$
\I_{P^{-1}\vf}\to I_\vf
$$

The element $\lambda$ defines a total order
on
$$
\I_\V:=\bigsqcup_\vf I_\vf.
$$

We have a natural map
$$
M:\I_\U=\bigsqcup\limits_{\vf}
\I_{P^{-1}\vf}\to \bigsqcup\limits_{\vf}
I_{\vf}=I_\V
$$

\begin{Lemma}
1) There is a unique total order
on $\I_\U $ such that:

--- the  map $M$ is monotonous;

--- this order agrees with those on each
$\I_{P^{-1}\vf}$;

\end{Lemma}

\begin{proof}
  If such an order exists, it must be defined as follows:

--- if $x,y\in \I_\U$, and $M(x)<M(y)$, then
$x<y$;

--- if $M(x)=M(y)\in I_\phi$, then $x,y\in
\I_{P^{-1}\vf}$.

It is clear that this way we indeed get a
total order on $\I_\U$. The map $M$ is
automatically monotonous. We only need to
check the matching of this order with that
on each $\I_{P^{-1}\vf}$. This follows
immediately from the monotonicity of the
corresponding maps
$$
\I_{P^{-1}\vf}\to I_\vf.
$$
\end{proof}

Next we define a map $W':\I_\U\to J$ as a
composition
$$
\I_\U\to \I_\V \stackrel{W}\to J
$$

\begin{Lemma} The constructed order on
$\I_\U$ and the constructed map $W'$ give
rise to an element in $\seq(\U)$ \end{Lemma}
\begin{proof}
  Straightforward
\end{proof}

We define the composition of the elements $\lambda_\vf$ and $\lambda$
to be the constructed element in $\seq(\U)$.

One can check that this composition
satisfies the associativity property.

\subsection{$\seq$-algebra structure on $\shom^\bullet$}
We need a couple of auxiliary constructions.
\subsubsection{} Given a dg- category $A$,
an ordinal $J$, and a map $X:J\to A$ we call
any element in $A(X)$ (see the very
beginning of the paper) {\em a chain in $A$}
or, more, specific, an $X$-chain in $A$. Fix
a chain $h\in A(X)$.

Suppose we are given an $X:J\to A$ as above.
Suppose that in addition, we are given an
ordinal $R$ and  an $R$-family of functors
$F_r:A\to B$.

Next,  for each $\vr=\vec{r_1r_2}\in \vR$,
choose  ordinals $I_\vr$ and elements
\begin{equation}\label{ash}
h_\vr\in \shom^{I_\vr}(F_{r_1};F_{r_2})
\end{equation}

Finally let us fix a monotonous map

$$W:I=\bigsqcup\limits_{\vr\in \vR} I_\vr\to J,$$
where the order on $I$ is defined by those
on $\vR$ and on each of $I_\vr$.

Given all these data, we will construct a
chain $c$ in $B$.

Before giving a formal definition let us
consider an example in which:

$$
J=\{0_j<1_j<\cdots<10_j\};
$$
$$
R=\{a<b<c\};
$$
$$
I_{\vec{ab}}=\{0_a<1_a<2_a<3_a\};\
I_{\vec{bc}}=\{0_b<1_b<2_b\},
$$

and the map $W$ is given by the following
table
$$
0_a\mapsto 1_j;\ 1_a,2_a\mapsto 3_j;\
3_a\mapsto 4_j;
$$
$$
0_b\mapsto 5_j;\ 1_b\mapsto 7_j;\ 2_b\mapsto
8_j.
$$

Then our map is constructed according to the
following picture:

$$
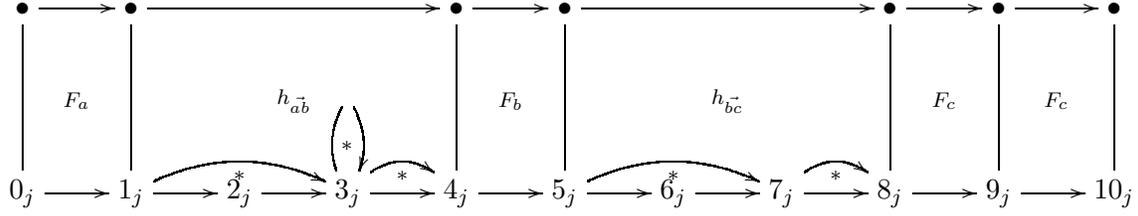
\xymatrix{\bullet\ar[r]&\bullet\ar[rrr]&&&\bullet\ar[r]&\bullet\ar[rrr]&&&\bullet\ar[r]&\bullet\ar[r]&\bullet\\
&&&\ar@/^0.5pc/[d]&&&&&&&\\
 0_j\ar[r]\ar@{-}[uu]\ar@{}[ruu]|{F_a} &
1_j\ar[r]\ar@/^1pc/[rr]_{*}
\ar@{-}[uu]\ar@{}[rrruu]|{h_{\vec{ab}}}&
2_j\ar[r] & 3_j\ar[r]
\ar@{-}@/^0.5pc/[u]_{*}\ar@/^1pc/[r]_{*} &
4_j\ar[r]\ar@{-}[uu]\ar@{}[ruu]|{F_b} &
5_j\ar[r]\ar@/^1pc/[rr]_{*}
\ar@{-}[uu]\ar@{}[rrruu]|{h_{\vec{bc}}}&
6_j\ar[r] & 7_j\ar[r]\ar@/^1pc/[r]_{*}&
8_j\ar[r]\ar@{-}[uu]\ar@{}[ruu]|{F_c} &
9_j\ar[r]\ar@{-}[uu]\ar@{}[ruu]|{F_c}
&10_j\ar@{-}[uu]}
$$

{\em  Explanation of the picture:} in each
cell marked by $*$ we compose the arrows on
the bottom of the cell;

in each cell marked by $h_{\vec{ab}}$,
$h_{\vec{bc}}$ we apply the corresponding
element from (\ref{ash});

in each cell marked by $F_a,F_b$ or $F_c$ we
apply the corresponding functor.

The resulting chain $c$ corresponds to the
top line of arrow on the picture.

\subsubsection{} Let us now make a formal
definition.

 First of all we need to
construct an ordinal $K$ and a map $Y:K\to
B$ so that $c\in B(Y)$.

{\em Constructing $K$}

 For $r\in R$ let $m_r$ be the supremum in
 $J$ of the set
 $$
\bigsqcup\limits_{\vec{r_1r_2}| r_2\leq r}
W(I_{\vec{r_1r_2}})
 $$

Let $M_r$ be the infimum in $J$ of the set
$$
\bigsqcup\limits_{\vec{r_1r_2}| r \leq r_1}
W(I_{\vec{r_1r_2}})
 $$

 We define:

 $$K:=K(J,W):=\bigsqcup\limits_{r\in R} [m_r,M_r].
$$
 We then have  natural maps
 $$
\pi:K\to R;
 $$
 $$
\kappa:K\to J.
 $$

{\em Constructing  a map $Y:=Y(X,W):K\to B$}
  Set
  $$
Y(j_r)=F_r(\kappa(j_r)),
  $$
where $j_r\in [m_r;M_r]\subset K$.

{\em Constructing the resulting chain $c\in
B(Y)$}

  We will define a map
  $$
\mu_W: A(X)\otimes
\bigotimes\limits_{\vec{r_1r_2}\in \vR}
\shom^{I_{\vec{r_1r_2}}}(F_{r_1};F_{r_2})\to
B(X).
  $$

so that $$c= \mu_W(h;\{h_\vr\}_{\vr\in
\vR}).
$$

  For an interval $[a,b]\subset J$, let
  $X_{ab}:=X|_{[a,b]}$.

Let $R=\{0<1<2<\cdots<N\}$.
  We then have
  $$
A(X)=A(X_{m_0M_0})\otimes
A(X_{M_0m_1})\otimes A(X_{m_1M_1})\otimes
\cdots .
  $$

Let $W_{r,r+1}:=W|_{I_{\vec{r,r+1}}}$.  We
then have a dominant map
$$
W_{r,r+1}:I_{\vec{r,r+1}}\to [M_r;m_{r+1}].
$$
Set
$$
X'_{r,r+1}:
I_{\vec{r,r+1}}\stackrel{W_{r,r+1}}\to
[M_r;m_{r+1}]\stackrel{X}\to A.
$$

Hence, we have an  induced map
$$
W_{r,r+1}^*: A(X_{M_r;m_{r+1}})\to
A(X'_{r,r+1});
$$
via substitution, we get a map:

$$
A(X_{M_rm_{r+1}})\otimes
\shom^{I_{\vec{r,r+1}}}(F_r;F_{r+1}) \to
B(F_r(X(M_r));F_{r+1}(X(m_{r+1})))
$$

The functor $F_r$ induces a map
$$
A(X_{m_rM_r})\to B(Y|_{[m_r;M_r]})
$$

Combining these maps we get the desired map
$ \mu_W. $

\subsubsection{Definition of a
$\seq$-algebra structure} Let us now
construct the structure maps
$$
A:\seq(\U)\to \full_{\shom^\bullet}(\U),
$$

Equivalently,  for each $\U$-diagram in
$\C$, one has to construct a map
$$
k[\seq(\U)]\otimes \bigotimes_{\vf\in \vF}
\shom^{I_\vf}(\D|_\vf)\to \shom^J(p_*\D),
$$
where  $\D|_\vf$ is the $\C$-globe obtained
by the restriction of $\D$ onto $\vf$ and
$p_*\D$ is a $\C$-globe obtained by
pre-composing $\D$ with the terminal map
$p:\U\to \G$.

Let $\mu$ be the minimal vector in $\vC$.
 Consider the ordinals $I_\vf$, $\vf \in
\vF_\mu$. It follows that they form
intervals in $\I$ and $\vf_1<\vf_2$ implies
$I_{\vf_1}<I_{\vf_2}$ in $\I$.

We have a restriction
$$
W:\bigsqcup\limits_{\vf\in \vF_\mu} I_\vf\to
J,
$$
satisfying all the conditions of the
previous section.  Hence we have a  map
$\mu_W$ as explained above.  Let us show
that after the application of $\mu_W$, the
remaining ingredients form a similar
structure to that with which we started.

The map $\mu_W$ only involves the complexes
$\shom^{I_\vf}$ with $\vf\in \vF_\mu$. The
remaining complexes are labelled by the
elements of the set
$$
\I':=\I\backslash \bigsqcup\limits_{\vf\in
\vF_\mu} I_\vf.
$$

The map $W$ naturally descends to a map
$$
W':\I'\to  K(W)
$$

Let $\C':=\C\backslash m_\C$;  We then get a
diagram $\U'$ with $\C_{\U'}=\C'$ and
$\F_{\U',\vc}=\F_{\U,\vc}$ It then follows
that $(\I',W')\in \seq(\U')$. Thus we have
constructed a map
$$
\nu_\U:\k[\seq(\U)]\otimes
\bigotimes\limits_{\vf\in \vF}
\shom^{I_\vf}(\D|_\vf)\to \k[\O(\U')]\otimes
\bigotimes\limits_{\vf'\in \vF_{\U'}}
\shom^{I_{\vf'}}(F_{P_*\D|_{\vf'}})
$$

One can now iterate this procedure thus
exhausting all the arguments; in the end we
will obtain a chain of morphisms in
$C_{M_\C}$, and, finally, we can take the
composition of all morphisms in this chain,
which will produce the result.

We omit the proof that this is indeed a
$\seq$-algebra structure --- this is pretty
clear.

\subsubsection{}

In Sec. \ref{consthom} we have defined a
map of co-simplicial complexes
$$
\hom^\bullet(F,G)\to \rhom^\bullet(F,G)
$$
for every pair of functors $F,G:A\to B$.
This way we get  a map of $\C$-complexes
$$
\hom^\bullet\to \rhom^\bullet.
$$

It is easy to check that $\hom^\bullet$ is a
$\seq$-subalgebra of $\rhom^\bullet$.
Furthermore, given a 2-ordinal $\U$ and an
$\U$-diagram $\D$,  for every $e\in
\seq(\U)$, the structure map
$$
\bigotimes_{\vf\in\vF_\U} \hom^{I_\vf}(\D_{\vf})
\xrightarrow{1\mapsto e} \seq(\U)\otimes\bigotimes_{\vf\in
\vF_\U}\hom^{I_\vf}(\D_{\vf}) \longrightarrow \hom^J(p_*\D)
$$
is the same.

This can be formulated as follows. Let $\T$
be the trivial $\nat$-colored 2-operad:  for
every $\nat$-colored 2-ordinal $\U$
$$
\T(\U):=\{1\},
$$
(the structure maps are then uniquely
defined); Let
\begin{equation}\label{trivproj}
 \seq \to \T
\end{equation}
be the obvious projection. We then have:

\begin{Proposition}\label{prophom}

1) $\hom^\bullet\subset \Rhom^\bullet$ is a
$\seq$-subalgebra;

2) the $\seq$-action on $\hom^\bullet$
passes through the projection
(\ref{trivproj}) \end{Proposition}

\subsubsection{Co-simplicial structure} Let
us recover the cosimplicial structure on
$\shom^\bullet$ from the $\seq$-structure.

Let $\G,I,J$ be the globe colored by the
ordinals $I$ and $J$. By definition,
$$
\seq(\G,I,J)=\Delta(I,J)
$$
is the space of all monotonous maps.

The 2-operadic structure gives rise to
associative maps
$$
\seq(\G,I,J)\times \seq(\G,J,K)\to
\seq(\G,I,K)
$$

 thus giving rise to a category structure on
 $\nat$ which is just given by composing the
 corresponding monotonous maps. That is,
 this category is nothing else but the
 simplicial category $\Delta$.
\subsubsection{}  Given a colored  2-ordinal $\U'$
with the underlying 2-ordinal $\U$ and a
coloring given by the family of ordinals
$I_\vf$, $\vf\in \vF_\U$; $J$, write

$$
\seq(\U)_{{\{I_\vf\}_{\vf\in
\vF}}}^J:=\seq(\U').
$$

As a function in $I_\vf,J$, $\seq(\U)$
becomes a functor
$$
(\Delta^\op)^{\vF}\times \Delta\to \sets.
$$

\subsection{Passing to the complexes}

 Define the realization
$$
\O(\U):=|\seq(\U)|:=\hom_{\Delta}(\S;
\k[\seq(\U)]\otimes_{(\Delta^\op)^\vF}
(\S)^{\boxtimes \vF}),
$$
where $\S$ is as in (\ref{simplchain}).

 It is immediate that these realizations
 form a dg-2-operad $\O$, and that this operad
 acts on the $\C$-complex
 $\Rhom=\hom_\k(\S,\shom^\bullet)$.

 Our goal now is to check that this operad
 satisfies the theorem.
\subsection{Contractibility of $\O$}

  First of all, let us
 construct a quasi-isomorphism
 $$
\O\to \triv
 $$

It is easy to construct such a map: it is
just the augmentation map. Let us show that
this map is a quasi-isomorphism

\begin{Proposition} For each ordinal $J$,
the poly-simplicial realization
$$
S(\U,J):=|\seq(\U)^J_{\bullet,\ldots,\bullet}|
$$
with respect to all lower indices is
contractible \end{Proposition}

\begin{proof}
  One can describe this realization explicitly. Let us so do: a point
  of $S(\U,J)$ is given by an equivalence class of the following data:

 1) a decomposition of  a fixed segment $I:=I_\U$ of length
 $|\vF|$ into a
 number of subsegments labelled by the
 elements from $\vF_\U$. The labelling should
 satisfy:

 a) if $\vf_1,\vf_2\in \vF$ and a segment
 labelled by $\vf_2$ lies between those
 labelled by $\vf_1$, then
 $\pi(\vf_1)>\pi(\vf_2)$;

 b) if $\pi(\vf_1)=\pi(\vf_2)=\vc$ and
 $\vf_1<\vf_2$ in $\vF_\vc$, then all
 segments labelled by $\vf_1$ are on the LHS
 of elements labelled by $\vf_2$;

 c) the total length of all segments
 labelled by the same element $\vf$ is 1

d) a monotonous map $\vJ\to I$.

Two such points are equivalent if one is
obtained from another by a number of
operations of the following two types:

--- adding into or deleting from our
decomposition a number of labelled segments
of length 0;

--- joining two neighboring segments of our
decomposition labelled by the same letter
into one segment labelled by the same
letter, or the inverse operation.

This space receives an obvious CW-structure.
The proof that this space is indeed a
realization is straightforward.

Let $S(\U)$ be the space whose points are
described by a)-c) (without d), and the
equivalence relation is the same. We then
get an obvious isomorphism
$$
S(\U,J)=S(\U)\times \Delta^J.
$$

{\bf Remarks}  1) One can prove that this is
an isomorphism of co-simplicial sets.

2) The topological realization
$$|\Sigma(\U)|:=|S(\U,\bullet)|$$ is then identified with the
space $S(\U)\times R$, where $R$ is the
space of monotonous maps of a unit segment
into the segment $I$. The spaces
$\Sigma(\U)$ form a  topological 2-operad.
This operad acts on  topological
realizations of $\shom^\bullet(F,G)$

\subsubsection{} Thus, it only remains to show
that the space $S(\U)$ is contractible.

For simplicity, let us identify $\C$ with
the set $0<1<\cdots <n$. Let $\U_m$ be a
ball in $\U$ which is the preimage of $
[m,n]\subset \C$. We then have a natural
projection
$$
P_m: S(\U_m)\to S(\U_{m+1})
$$

this projection sends a point in $S(\U_m)$
into a point obtained from it by collapsing
each segment labelled by elements from
$\pi^{-1}\vec{m,m+1}$ to a point.

Conversely, given:

-- a point $x\in S(\U_{m+1})$;

-- a monotonous map
$U:\pi^{-1}(\vec{m,m+1})\to I_{\U_{m+1}}$,

one can construct a point in $S(\U_m)$ by
inserting unit segments  labelled by $i\in
\pi^{-1}(\vec{m,m+1})$ in place of the point
$U(i)$.

It is clear that this way we get a bijection
$$
S(\U_m)\cong S(t_{m+1})\times
\Delta^{\pi^{-1}m}.
$$

This argument implies that the space $S(\U)$
is homeomorphic to a product of simplices,
hence is contractible.

\end{proof}

\subsubsection{}  Thus, we have proven the
assertion of the Theorem \ref{maintheor}
that a contractible 2-operad acts on
$\Rhom$. It remains to check the conditions
1,2. They follow immediately from
Proposition \ref{prophom}. This completes
the proof of the Theorem.

\section{Relation to Deligne's conjecture}

Given a dg-category $\A$, we can consider a
complex $\rhom(\Id_\A,\Id_\A)$ This complex
is called the Hochschild complex of the
category $\A$. If the category $\A$ has only
one object $p$, then its Hochschild complex
coincides with that of the associative
algebra $\End_\A(p)$.

Thus, we denote $$
\hoch_\A:=\rhom(\Id_\A,\Id_\A).$$

As a corollary of  the just proven theorem,
we have a certain structure on $\hoch_\A$;
before defining it, let us give it a name
"an $\O$-algebra structure on $\hoch_\A$"

The definition is as follows. Given a
complex $\K$ (for example $\K:=\hoch_\A$) we
define a 2-operad $\full_\K$ by setting
$$
\full_\K(\U):=\hom_\k(\K^{\otimes \F_\U};\K)
$$
with the obvious insertion maps.

{\bf Remark} Of course, this construction is
a particular case of the full 2-operad of a
$\C$-complex, where $\C$ is the category
with one object and one arrow so that there
is only one globe in $\C$, and a
$\C$-complex is the same as a usual complex,
so that our complex $\K$ gives rise to a
$\C$-complex and we can apply the
construction of the full operad of a
$\C$-complex. This way we get another
construction of $\full_\K$.

Given a 2-operad $\E$, we define {\em an
$\E$-algebra structure on $\K$} as a map of
2-operads
$$
E\to \full_\K.
$$

Theorem \ref{maintheor} immediately implies
that:
\begin{Proposition}
 $\hoch_\A$ has a structure of algebra
over the  2-operad $\O$, as in the statement
of the Theorem.
\end{Proposition}

As $\O$ is a contractible 2-operad, a result
from \cite{Bat} readily implies that
\begin{Corollary}
 A
certain operad which is homotopy equivalent
to the chain operad of little disks acts on
$\hoch_\A$. \end{Corollary}

This corollary is known as {\em Deligne's
conjecture on Hochschild cochains}.

\end{document}